\documentclass[11pt]{article}
\usepackage{amsfonts,latexsym,rawfonts,amsmath,amssymb,amsthm,graphicx}
\numberwithin{equation}{section}
\newtheorem{theorem}{Theorem}[section]
\newtheorem{lem}[theorem]{Lemma}
\newtheorem{thm}[theorem]{Theorem}

\newtheorem{rem}[theorem]{Remark}

\def\endproof{$\hfill\Box$\\}

\title{ $L^p$ Liouville type theorems for harmonic functions on gradient Ricci solitons}
\author{Yong Luo}
\date{}
\begin{document}
\maketitle
\begin{abstract}
In this paper we consider $L^p$ Liouville type theorems for harmonic functions on gradient Ricci solitons. In particular, assume that $(M,g)$ is a gradient shrinking or steady  K\"ahler-Ricci soliton, then we prove that any pluriharmonic function $u$ on $M$ with $\nabla u\in L^p(M)$ for some $1<p\leq 2$  is a constant function.
\end{abstract}
\section{Introduction}
The classical Liouville theorem states that any bounded holomorphic function on the complex plane $\mathbb{C}$ must be constant, which also holds for harmonic functions on $\mathbb{R}^n$. In 1975, Yau \cite{Yau1} generalized this theorem by proving that any positive harmonic function on an $n$-dimension complete manifold with nonnegative Ricci curvature must be constant, and later Cheng and Yau \cite{ChY} showed a gradient estimate on the manifold with Ricci curvature lower bound to give an effective version of the Liouville theorem. Cheng and Yau's gradient estimate method plays a crucial role in the study of many partial differential equations on manifolds (see e.g. \cite{CM, Li} and references therein for more information).

 It is also very interesting to consider Liouville type theorem under integrability conditions. For example, it's natural to ask whether a manifold possesses a non-trivial $L^2$ harmonic function. In 1976, Yau \cite{Yau2} proved that any $L^p(p>1)$ harmonic function on a complete manifold must be constant. As a corollary, by using the Bochner formula, he proved that any harmonic function $u$ with $\nabla u$ in $L^p(p>1)$ on a complete manifold with nonnegative Ricci curvature must be constant.

 In this paper we are curious with Liouville type theorems for harmonic functions on gradient Ricci solitons. A complete Riemannian manifold $(M^n,g)$ is a gradient Ricci soliton if there exists a smooth function $f$ on $M$ satisfying the equation
 \begin{equation}
\rm Ric+\nabla\nabla f=\mu g,
 \end{equation}
where $\mu$ is $\frac12$, or $0$, or $-\frac12$. It is called a gradient shrinking Ricci soliton if $\mu=\frac12$, or a gradient steady Ricci soliton if $\mu=0$, or a gradient expanding Ricci soliton if $\mu=-\frac12$. Gradient Ricci solitons could be seem as natural generalizations of Einstein manifolds, since when $f$ is a constant function, the gradient Ricci soliton is simply an Einstein manifold. On the other hand, gradient Ricci solitons are also self-similar solutions to Hamilton's Ricci flow and play a crucial role in the study of formation of singularities in the Ricci flow (see \cite{Cao}).

 It is well known that gradient shrinking or steady Ricci solitons share many properties with complete manifolds with nonnegative Ricci curvature. Therefore it is very natural to ask whether the above mentioned Liouville type theorems  for harmonic functions due to Yau on complete manifold with nonnegative Ricci curvature still hold on shrinking or steady Ricci solitons. In fact much effort has been made along this direction in the past decade.

Munteanu and Wang \cite{MW2} showed that any bounded holomorphic function on shrinking K\"ahler-Ricci  solitons must be constant. Mai and Ou \cite{MO} proved that any bounded harmonic function on gradient shrinking  Ricci solitons with constant scalar curvature must be constant. 
Note that though we already have these advances, the following problem remains a challenge.
\vspace{0.2cm}
\\$\textbf{Problem 1}$. Is any positive harmonic function on gradient shrinking Ricci  solitons a constant function?

\vspace{0.2cm}

There were also several results on $L^p$ type Liouville theorems for harmonic functions on gradient shrinking or steady  Ricci solitons. Munteanu and Sesum \cite{MS} proved that a harmonic function $u$ on a gradient  shrinking K\"ahler-Ricci soliton with $\nabla u\in L^2$ is a constant function. In the same paper they proved that a harmonic function $u$ on a gradient  steady  Ricci soliton with $\nabla u\in L^2$ is a constant function. Dung, Le Hai and Thanh \cite{DHT} proved that a harmonic function $u$ on a gradient steady  Ricci soliton $(M^n,g)$ with $\nabla u\in L^p(2<p<\frac{4n}{2n-1})$ is a constant function, if $(M,g)$ satisfies the Euclidean volume growth condition. Inspired by Yau's $L^p$ Liouville type theorem for harmonic functions on manifolds with nonnegative Ricci curvature and these results, it is very natural to study the following problem.
\vspace{0.2cm}
\\$\textbf{Problem 2}$. Is any harmonic function $u$ with $\nabla u\in L^p(p>1)$ on gradient  shrinking or steady Ricci solitons a constant function?
\vspace{0.2cm}

For gradient shrinking K\"ahler-Ricci  solitons we answer Problem 2 partially in the affirmative, which extends Munteanu and Sesum's result.
\begin{theorem}\label{main thm1}
Let $(M,g)$ be a gradient shrinking K\"ahler-Ricci soliton. If $u$ is a pluriharmonic function with $\int_M|\nabla u|^p<\infty(1<p\leq2)$ , then $u$ is a constant function.
\end{theorem}
For gradient steady K\"ahler-Ricci solitons  we answer Problem 2 partially in the affirmative too.
\begin{theorem}\label{main thm2}
Let $(M,g)$ be a gradient steady K\"ahler-Ricci soliton. If $u$ is a pluriharmonic function with $\int_M|\nabla u|^p<\infty(1<p\leq 2)$ , then $u$ is a constant function.
\end{theorem}
We remark that Liouville type theorems for $f$-harmonic functions on gradient shrinking Ricci solitons were studied by Ge and Zhang \cite{GZ} and Wu and Wu \cite{WW} and others, which are relatively simpler and complete.

At last, we have the following result.
\begin{thm}\label{application}
Let $(M,g)$ be a gradient  K\"ahler-Ricci soliton with potential function $f$. Then if $f$ is a harmonic function on $M$ with $\int_M|\nabla f|^p<\infty(p>1)$, $(M,g)$ is a K\"ahler Einstein manifold.
\end{thm}


\vspace{0.2cm}
 \textbf{Organization.} Theorem \ref{main thm1} and Theorem \ref{application} are proved in section 2, and  Theorem \ref{main thm2} is proved in section 3. In the appendix we discuss more $L^p$ Liouville type theorems for harmonic or pluriharmonic functions on gradient shrinking or expanding Ricci solitons.

\section{Gradient K\"ahler-Ricci solitons}
Let $(M,g)$ be a  complete noncompact  gradient  K\"ahler-Ricci soliton, that is
$$R_{\alpha\bar{\beta}}+f_{\alpha\bar{\beta}}=\mu g_{\alpha\bar{\beta}},$$
for a smooth function $f$ on $M$ satisfying
\begin{eqnarray}\label{eq1}
f_{\alpha\beta}=f_{\bar{\alpha}\bar{\beta}}=0.
\end{eqnarray}
Here and in other places of this paper we use notations
$$u_\alpha:=\frac{\partial u}{\partial z_{\alpha}},\ u_{\alpha\beta}=\nabla_{\partial z_\beta}\frac{\partial u}{\partial z_\alpha}, \ u_{\bar{\alpha}}=\frac{\partial u}{\partial \bar{z}_\alpha}, \cdot\cdot\cdot.$$

When $\mu=1$, it is well known that after normalizing the potential function we have ( see \cite{CLN}\cite{ELM})
\begin{eqnarray}\label{eq0}
|\nabla f|^2+R=f,
\end{eqnarray}
and Chen \cite{Chen} proved that $R\geq0$.

For nontrivial gradient shrinking Ricci solitons, From \cite{CZ} and \cite{HM} we have
\begin{eqnarray}\label{ine0}
\frac14(r(x)-c)^2\leq f(x)\leq\frac14(r(x)+c)^2.
\end{eqnarray}
Here and in the rest of this paper we will use $C$ or $c$ to denote universal constants, which may be different with each other.

From (\ref{eq0}) and (\ref{ine0}) we have
\begin{eqnarray}\label{ine00}
|\nabla f|\leq\frac12(r(x)+c).
\end{eqnarray}

Munteanu and Sesum \cite{MS} proved that any harmonic function on $(M,g)$ with bounded total energy is a constant function. We will extend their result to the following theorem.
\begin{thm}\label{thm2.1}
Let $(M,g)$ be a gradient shrinking K\"ahler-Ricci soliton. If $u$ is a pluriharmonic function on $M$ with $\int_M|\nabla u|^p<\infty$ for some constant $1<p\leq 2$ , then $u$ is a constant function.
\end{thm}
In their proofs, Munteanu and Sesum used  a result due to Li \cite{Li2} that if $u$ is a harmonic function and with finite total energy on a K\"ahler manifold then it is in fact pluriharmonic. We extend Li's result to the following, which will deduce Theorem \ref{application}.
\begin{lem}\label{main lem1}
Assume that $u$ is a harmonic function on a complete K\"ahler manifold $(M,g)$ with  $\int_M|\nabla u|^p<\infty$ for some constant $p>1$ and $u_{\alpha\beta}$=0 for any $\alpha,\ \beta$, then $u$ is a pluriharmonic function.
\end{lem}
\proof Let $\phi: M\to [0,1]$ be a cut-off function such that $\phi=1$ on $B_p(r)$(a geodesic ball centered at $p\in M$ of radius $r$) and $|\nabla\phi|\leq\frac{C}{r}$. Assume that $\alpha$ is constant to be determined later and $\epsilon$ is a positive number. Note that
\begin{eqnarray*}
|\nabla u|^2&=&u_{\alpha}u_{\bar{\alpha}},
\\ \Delta u&=&u_{\alpha\bar{\alpha}}.
\end{eqnarray*}
Then
\begin{eqnarray}\label{ine1}
&&\int_M\phi^2(|\nabla u|^2+\epsilon)^\alpha|u_{\alpha\bar{\beta}}|^2\nonumber
\\&=&-2\int_M\phi\phi_{\bar{\beta}}(|\nabla u|^2+\epsilon)^\alpha u_{\alpha}u_{\bar{\alpha}\beta}-\int_M\phi^2(|\nabla u|^2+\epsilon)^\alpha u_\alpha u_{\bar{\alpha}\beta\bar{\beta}}\nonumber
\\&-&\alpha\int_M\phi^2(|\nabla u|^2+\epsilon)^{\alpha-1}u_{\alpha}u_{\bar{\alpha}\beta}[u_{\gamma\bar{\beta}}u_{\bar{\gamma}}+u_{\gamma}u_{\bar{\gamma}\bar{\beta}}]\nonumber
\\&=&-2\int_M\phi\phi_{\bar{\beta}}(|\nabla u|^2+\epsilon)^\alpha u_{\alpha}u_{\bar{\alpha}\beta}-\int_M\phi^2(|\nabla u|^2+\epsilon)^\alpha u_\alpha u_{\beta\bar{\beta}\bar{\alpha}}\nonumber
\\&-&\alpha\int_M\phi^2(|\nabla u|^2+\epsilon)^{\alpha-1}u_{\alpha}u_{\bar{\alpha}\beta}u_{\gamma\bar{\beta}}u_{\bar{\gamma}}\nonumber
\\&\leq&\frac1\delta\int_M(|\nabla u|^2+\epsilon)^\alpha|\nabla u|^2|\nabla\phi|^2+\delta\int_M\phi^2(|\nabla u|^2+\epsilon)^\alpha|u_{\alpha\bar{\beta}}|^2\nonumber
\\&-&\alpha\int_M\phi^2(|\nabla u|^2+\epsilon)^{\alpha-1}u_{\alpha}u_{\bar{\alpha}\beta}u_{\gamma\bar{\beta}}u_{\bar{\gamma}},
\end{eqnarray}
where in the second equality we used $$u_{\bar{\alpha}\beta\bar{\beta}}=u_{\beta\bar{\beta}\bar{\alpha}}$$ since $(M,g)$ is a K\"ahler manifold. Note that 
\begin{eqnarray*}
0\leq u_{\alpha}u_{\bar{\alpha}\beta}u_{\gamma\bar{\beta}}u_{\bar{\gamma}}=|u_{\alpha}u_{\bar{\alpha}\beta}|^2.
\end{eqnarray*}

If $p\geq2$, we let $\alpha=\frac{p-2}{2}\geq0$ and $\delta=\frac12$.
we get
\begin{eqnarray*}
&&\int_M\phi^2(|\nabla u|^2+\epsilon)^\alpha|u_{\alpha\bar{\beta}}|^2
\\&\leq&2\int_M(|\nabla u|^2+\epsilon)^\alpha|\nabla u|^2|\nabla\phi|^2+\frac12\int_M\phi^2(|\nabla u|^2+\epsilon)^\alpha|u_{\alpha\bar{\beta}}|^2,
\end{eqnarray*}
which implies that
\begin{eqnarray*}
&&\int_M\phi^2(|\nabla u|^2+\epsilon)^\alpha|u_{\alpha\bar{\beta}}|^2
\\&\leq&4\int_M(|\nabla u|^2+\epsilon)^\alpha|\nabla u|^2|\nabla\phi|^2,
\end{eqnarray*}
and hence
\begin{eqnarray*}
&&\int_M\phi^2|\nabla u|^{p-2}|u_{\alpha\bar{\beta}}|^2
\\&\leq&4\int_M|\nabla u|^p|\nabla\phi|^2.
\end{eqnarray*}
Letting $r\to \infty$ we get $u_{\alpha\bar{\beta}}=0$, i.e. $u$ is a pluriharmonic function.

If $1<p<2$, let $\alpha<0$ to be determined later, then from (\ref{ine1}) we have
\begin{eqnarray*}
&&\int_M\phi^2(|\nabla u|^2+\epsilon)^\alpha|u_{\alpha\bar{\beta}}|^2
\\&\leq&\frac1\delta\int_M(|\nabla u|^2+\epsilon)^\alpha|\nabla u|^2|\nabla\phi|^2+\delta\int_M\phi^2(|\nabla u|^2+\epsilon)^\alpha|u_{\alpha\bar{\beta}}|^2\nonumber
\\&-&\alpha\int_M\phi^2(|\nabla u|^2+\epsilon)^{\alpha-1}|\nabla u|^2|u_{\alpha\bar{\beta}}|^2.
\end{eqnarray*}
Then by letting $\alpha=\frac{p-2}{2}$ and $\delta=\frac{2-p}{2}$ that
\begin{eqnarray*}
&&(p-1)\int_M\phi^2(|\nabla u|^2+\epsilon)^\alpha|u_{\alpha\bar{\beta}}|^2
\\&\leq&\frac{2}{(2-p)}\int_M(|\nabla u|^2+\epsilon)^\alpha|\nabla u|^2|\nabla\phi|^2\nonumber
\\&\leq&\frac{2}{(2-p)}\int_M|\nabla u|^p|\nabla\phi|^2.
\end{eqnarray*}
Letting $r\to \infty$ we get $u_{\alpha\bar{\beta}}=0$, i.e. $u$ is a pluriharmonic function. \endproof

Theorem \ref{application} can be deduced as follows. By the above lemma, $f$ is a pluriharmonic function, therefore we must have $$R_{\alpha\bar{\beta}}=\mu g_{\alpha\bar{\beta}},$$ i,e. $(M,g)$ is a K\"ahler Eistein manifold.

Next we are going to prove Theorem \ref{thm2.1} by following Munteanu and Sesum's argument with several modifications.

\proof We only need to prove case of $1<p<2$. Define
$$F:=\langle\nabla f,\nabla u\rangle=\frac12(u_\alpha f_{\bar{\alpha}}+u_{\bar{\alpha}}f_{\alpha}).$$
We will show that $F=0$. By (\ref{eq1}) and Lemma 2.2, we have
$$(u_{\alpha} f_{\bar{\alpha}})_{\bar{\delta}}=u_{\alpha\bar{\delta}}+u_{\alpha} f_{\bar{\alpha}\bar{\delta}}=0.$$
Similarly we have $$(u_{\bar{\alpha}}f_{\alpha})_{\delta}=0.$$ Therefore
$\Delta F=0$,
and
\begin{eqnarray*}
 \Delta F^2=2|\nabla F|^2.
\end{eqnarray*}
Let $v=F^2$, then
\begin{eqnarray}\label{ine2}
\Delta v= 2|\nabla F|^2.
\end{eqnarray}
 In the rest of this paper we will let $\phi$ be a cut-off function on $M$ as in the proof of the Lemma 2.2. Let $\alpha<0$ be a constant to be determined later. Then by (\ref{ine2}) and integration by parts we have
\begin{eqnarray*}
\int_M\phi^2(F^2+\epsilon)^\alpha|\nabla F|^2&=&\frac12\int_M\phi^2(v+\epsilon)^\alpha\Delta v
\\&=&-\int_M\phi(v+\epsilon)^\alpha\nabla\phi\nabla v-\frac12\alpha\int_M(v+\epsilon)^{\alpha-1}|\nabla v|^2
\\&=&-2\int_M\phi F(F^2+\epsilon)^\alpha\nabla\phi\nabla F-2\alpha\int_M(F^2+\epsilon)^{\alpha-1} F^2|\nabla F|^2
\\&\leq&-2\int_M\phi F(F^2+\epsilon)^\alpha\nabla\phi\nabla F-2\alpha\int_M(F^2+\epsilon)^\alpha|\nabla F|^2.
\end{eqnarray*}
Therefore
\begin{eqnarray*}
&&(1+2\alpha)\int_M\phi^2(F^2+\epsilon)^\alpha|\nabla F|^2
\\&\leq&\delta\int_M\phi^2(F^2+\epsilon)^\alpha|\nabla F|^2+\frac{1}{\delta}\int_MF^2(F^2+\epsilon)^\alpha|\nabla\phi|^2.
\end{eqnarray*}
We choose $\alpha=\frac{p-2}{2}$, and $\delta=\frac{p-1}{2}$, then
\begin{eqnarray*}
\frac{p-1}{2}\int_M\phi^2(F^2+\epsilon)^\alpha|\nabla F|^2\leq\frac{2}{p-1}\int_MF^2(F^2+\epsilon)^\alpha|\nabla\phi|^2.
\end{eqnarray*}
We obtain
\begin{eqnarray}\label{ine3}
\int_M\phi^2(F^2+\epsilon)^{\frac{p-2}{2}}|\nabla F|^2&\leq&(\frac{2}{p-1})^2\int_M|F|^p|\nabla\phi|^2
\\&\leq&\frac{C}{r^{2-p}}\int_M|\nabla u|^p,\nonumber
\end{eqnarray}
where we used $$|F|\leq|\nabla f||\nabla u|\leq Cr|\nabla u|$$ when $r$ is large enough, since $f$ satisfies (\ref{ine00}) . Then by letting $r\to \infty$ we have $\nabla F=0$, which yields that $F$ is a constant function on $M$.

Since $f$ growths to infinity when $r\to\infty$ by (\ref{ine0}), $f$ must attains its minimum at some point on $M$. Therefore $c=0$, i.e. $F\equiv 0$ on $M$. Define the $f$-Laplacian of a function to be $$\Delta_f=\Delta-\nabla f \cdot\nabla.$$
We now prove the following:
\begin{eqnarray}\label{f-subharmonic}
\Delta_f|\nabla u|\geq\frac12|\nabla u|,
\end{eqnarray}
by using Munteanu and Sesum's argument. Since we have proved $\langle\nabla f,\nabla u\rangle=0$, it follows that
$$\Delta_fu=\Delta u-\langle\nabla f,\nabla u\rangle=0.$$
Then the Bochner formula implies that
\begin{eqnarray*}
\Delta_f|\nabla u|^2&=&2Ric_f(\nabla u,\nabla u)+2\langle\nabla\Delta_fu, \nabla u\rangle+2|\nabla\nabla u|^2
\\&=&|\nabla u|^2+2|\nabla\nabla u|^2\geq|\nabla u|^2+2|\nabla|\nabla u||^2,
\end{eqnarray*}
where in the last inequality we used Kato inequality. Since on the other hand
$$\Delta_f|\nabla u|^2=2|\nabla u|\Delta_f|\nabla u|+2|\nabla|\nabla u||^2,$$
it is clear that we get (\ref{f-subharmonic}). Therefore $\Delta_f|\nabla u|\geq0$. Note that by the  assumption that $$\int_M|\nabla u|^p<\infty$$ and (\ref{ine0}) we have $$\int_M|\nabla u|^pe^{-f}<\infty.$$ Then from Theorem 1.1 in \cite{PRS} we can deduce that $u$ is a constant function. We also can give a direct proof as follows.

Note that
$$\Delta_f|\nabla u|^2=2\Delta_f|\nabla u||\nabla u|+2|\nabla|\nabla u||^2\geq 2|\nabla|\nabla u||^2.$$

Define $w=|\nabla u|^2$ and let $\alpha<0$ be a constant to be determined later, then
\begin{eqnarray*}
&&\int_M\phi^2(w+\epsilon)^\alpha|\nabla|\nabla u||^2e^{-f}
\\&\leq&\frac12\int_M\phi^2(w+\epsilon)^\alpha\Delta_fwe^{-f}
\\&=&-\int_M\phi\nabla\phi\nabla w(w+\epsilon)^\alpha e^{-f}-\frac12\alpha\int_M\phi^2(w+\epsilon)^{\alpha-1}|\nabla w|^2e^{-f}
\\&=&-\int_M\phi\nabla\phi |\nabla u|\nabla|\nabla u|(w+\epsilon)^\alpha e^{-f}-2\alpha\int_M\phi^2w(w+\epsilon)^{\alpha-1}|\nabla|\nabla u||^2e^{-f}
\\&\leq&-\int_M\phi\nabla\phi |\nabla u|\nabla|\nabla u|(w+\epsilon)^\alpha e^{-f}-2\alpha\int_M\phi^2(w+\epsilon)^{\alpha}|\nabla|\nabla u||^2e^{-f},
\end{eqnarray*}
Let $\alpha=\frac{p-2}{2}$, we see that
\begin{eqnarray*}
(p-1)\int_M\phi^2(w+\epsilon)^\alpha|\nabla|\nabla u||^2e^{-f}&\leq&\frac{p-1}{2}\int_M\phi^2(w+\epsilon)^\alpha|\nabla|\nabla u||^2e^{-f}
\\&+&\frac{2}{p-1}\int_M|\nabla u|^2(w+\epsilon)^\alpha|\nabla\phi|^2e^{-f}
\\&\leq&\frac{p-1}{2}\int_M\phi^2(w+\epsilon)^\alpha|\nabla|\nabla u||^2e^{-f}
\\&+&\frac{2}{p-1}\int_M|\nabla u|^p|\nabla\phi|^2e^{-f}.
\end{eqnarray*}
Therefore
\begin{eqnarray}
\int_M\phi^2(w+\epsilon)^\alpha|\nabla|\nabla u||^2e^{-f}\leq\frac{4}{(p-1)^2}\int_M|\nabla u|^p|\nabla\phi|^2e^{-f}.
\end{eqnarray}
Letting $r\to \infty$ in the above inequality we get $|\nabla|\nabla u||=0$, and therefore $|\nabla u|=c$. Since $M$ has infinite volume we see that $c=0$, and $u$ is a constant function.
\endproof
\begin{rem}\label{rem1}
We can see that in the proof of Theorem \ref{thm2.1}, the shrinking Ricci soliton equation was not essentially used. Assume that $Vol(M)=\infty. $To prove $F=0$, we only need $f_{\alpha\beta}=f_{\bar{\alpha}\bar{\beta}}=0$, $|\nabla f|$ growth linearly and $\nabla f$ has an interior zero point. While at the last step, we used a $L^p(p>1)$ Liouville type theorem for $\Delta_f|\nabla u|\geq0$, where we need deduce that $\int_M|\nabla u|^pe^{-f}<\infty$ from the assumption $\int_M|\nabla u|^p<\infty$, which holds if $f$ has a lower bound.
\end{rem}
From Remark \ref{rem1} we have the following $L^p$ Liouville type theorem for pluriharmonic functions on smooth K\"ahler measure metric spaces.
\begin{thm}
Let $(M,g,f)$ be a smooth K\"ahler measure metric space such that $Vol(M)=\infty$, $Ric_f=Ric+\nabla\nabla f\geq0$ and  $f$ satisfies the following conditions:
\\ \rm (a) $f_{\alpha\beta}=f_{\bar{\alpha}\bar{\beta}}=0$;
\\ \rm (b) $|\nabla f(x)|\leq Cr(x)+o(r(x)) $ when $r(x)$ is large enough, \ and \ $\nabla f$ has an interior zero point;
\\ \rm (c) $f$ is bounded from below.

Then if $u$ is a pluriharmonic function on $M$ with $\int_M|\nabla u|^p<\infty$ for some constant $1<p\leq 2$ , $u$ is a constant function.
\end{thm}
\section{Gradient steady  K\"ahler-Ricci solitons}
Let $(M,g)$ be a complete noncompact gradient steady K\"ahler-Ricci soliton, that is
$$R_{\alpha\bar{\beta}}+f_{\alpha\bar{\beta}}=0,$$
for a smooth function $f$ on $M$ satisfying
$$f_{\alpha\beta}=f_{\bar{\alpha}\bar{\beta}}=0.$$
Hamilton \cite{Ham} showed that $$R+|\nabla f|^2=\lambda,$$
for some constant $\lambda>0$. Since $R\geq0$(see \cite{Chen}), we see that $|\nabla f|\leq\sqrt{\lambda}$. In this section we will prove the following theorem.
\begin{thm}\label{main thm steady}
Let $(M,g)$ be a gradient steady K\"ahler-Ricci soliton, and assume that $u$ is a pluriharmonic function on $M$ with $\int_M|\nabla u|^p<\infty$ for some constant $1<p\leq 2$, then $u$ is a constant function.
\end{thm}

\proof We only need to prove case of $1<p<2$. Define
$$F:=\langle\nabla f,\nabla u\rangle=\frac12(u_\alpha f_{\bar{\alpha}}+u_{\bar{\alpha}}f_{\alpha}).$$
We show that $F=0$. Since $u$ is a pluriharmonic function, we can see that $\Delta F=0$, i.e. $F$ is a harmonic function on $M$. Then similarly with the proof of (\ref{ine3}) we can obtain
\begin{eqnarray*}
\int_M\phi^2(F^2+\epsilon)^{\frac{p-2}{2}}|\nabla F|^2&\leq&(\frac{2}{p-1})^2\int_M|F|^p|\nabla\phi|^2
\\&\leq&\frac{C}{r^2}\int_{B_p(2r)\setminus B_p(r)}|\nabla u|^p\nonumber
\end{eqnarray*}
where we used $$|F|\leq|\nabla f||\nabla u|\leq \sqrt{\lambda}|\nabla u|.$$ Then by letting $r\to \infty$ we have $\nabla F=0$, which yields that $F$ is a constant function on $M$. Assume that $F=c$. Then by
\begin{eqnarray}
|c|=|F|\leq|\nabla f||\nabla u|\leq\sqrt{\lambda}|\nabla u|,
\end{eqnarray}
we get $|\nabla u|\geq\frac{|c|}{\sqrt{\lambda}}$. Therefore $c=0$, since $\int_M|\nabla u|^p<\infty$ and $M$ has infinite volume \cite{MS}.

Let $\alpha<0$ be a constant to be determined later. Then
\begin{eqnarray}\label{eq2}
&&\int_M Ric(\nabla u,\nabla u)\phi^2(|\nabla u|^2+\epsilon)^\alpha\nonumber
\\&=&\int_MR_{ij}u_iu_j\phi^2(|\nabla u|^2+\epsilon)^\alpha
=-\int_Mf_{ij}u_iu_j\phi^2(|\nabla u|^2+\epsilon)^\alpha\nonumber
\\&=&\int_Mf_iu_{ij}u_j\phi^2(|\nabla u|^2+\epsilon)^\alpha+\int_Mf_iu_iu_j(\phi^2)_j(|\nabla u|^2+\epsilon)^\alpha\nonumber
\\&+&\alpha\int_Mf_iu_iu_j\phi^2(|\nabla u|^2+\epsilon)^{\alpha-1}(|\nabla u|^2)_j\nonumber
\\&=&\int_Mf_iu_{ij}u_j\phi^2(|\nabla u|^2+\epsilon)^\alpha,
\end{eqnarray}
where in the last equality we used $$f_iu_i=\langle\nabla f,\nabla u \rangle=F=0.$$
Note that
\begin{eqnarray*}
&&\int_Mf_iu_{ij}u_j\phi^2(|\nabla u|^2+\epsilon)^\alpha
\\&=&-\int_M\Delta fu_ju_j\phi^2(|\nabla u|^2+\epsilon)^\alpha-\int_Mf_iu_ju_{ij}\phi^2(|\nabla u|^2+\epsilon)^\alpha
\\&-&\int_Mf_iu_ju_j(\phi^2)_i(|\nabla u|^2+\epsilon)^\alpha-2\alpha\int_Mf_iu_ju_j\phi^2(|\nabla u|^2+\epsilon)^{\alpha-1}u_{ki}u_k.
\end{eqnarray*}
Then by $\Delta f=-R$, we get
\begin{eqnarray}\label{eq4}
&&\int_Mf_iu_{ij}u_j\phi^2(|\nabla u|^2+\epsilon)^\alpha\nonumber
\\&=&\frac12\int_MRu_ju_j\phi^2(|\nabla u|^2+\epsilon)^\alpha \nonumber
\\&-&\frac12\int_Mf_iu_ju_j(\phi^2)_i(|\nabla u|^2+\epsilon)^\alpha-\alpha\int_Mf_iu_ju_j\phi^2(|\nabla u|^2+\epsilon)^{\alpha-1}u_{ki}u_k.
\end{eqnarray}
Plug this in formula (\ref{eq2}) we obtain that
\begin{eqnarray}\label{eq3}
&&\int_M Ric(\nabla u,\nabla u)\phi^2(|\nabla u|^2+\epsilon)^\alpha\nonumber
\\&=&\frac12\int_MRu_ju_j\phi^2(|\nabla u|^2+\epsilon)^\alpha\nonumber
\\&-&\frac12\int_Mf_iu_ju_j(\phi^2)_i(|\nabla u|^2+\epsilon)^\alpha-\alpha\int_Mf_iu_ju_j\phi^2(|\nabla u|^2+\epsilon)^{\alpha-1}u_{ki}u_k.
\end{eqnarray}
Recall that we have the following Bochner formula
\begin{eqnarray*}
&&\Delta|\nabla u|^2\\&=&2Ric(\nabla u,\nabla u)+2|u_{ij}|^2
\\&\geq&2Ric(\nabla u,\nabla u)+2|\nabla|\nabla u||^2,
\end{eqnarray*}
where we used the Kato inequality
$$|u_{ij}|^2\geq |\nabla|\nabla u||^2.$$
Then we have
\begin{eqnarray*}\label{ine4}
&&\int_M\Delta|\nabla u|^2\phi^2(|\nabla u|^2+\epsilon)^\alpha
\\&\geq&2\int_MRic(\nabla u,\nabla u)\phi^2(|\nabla u|^2+\epsilon)^\alpha+2\int_M|\nabla|\nabla u||^2\phi^2(|\nabla u|^2+\epsilon)^\alpha.
\end{eqnarray*}
Plug this inequality into (\ref{eq3}) we get
\begin{eqnarray*}
&&2\int_M|\nabla|\nabla u||^2\phi^2(|\nabla u|^2+\epsilon)^\alpha+\int_MRu_ju_j\phi^2(|\nabla u|^2+\epsilon)^\alpha\nonumber
\\&\leq&\int_Mf_iu_ju_j(\phi^2)_i(|\nabla u|^2+\epsilon)^\alpha+2\alpha\int_Mf_iu_ju_j\phi^2(|\nabla u|^2+\epsilon)^{\alpha-1}u_{ki}u_k
\\&+&\int_M\Delta|\nabla u|^2\phi^2(|\nabla u|^2+\epsilon)^\alpha
\\&\leq&\int_Mf_iu_ju_j(\phi^2)_i(|\nabla u|^2+\epsilon)^\alpha+2\alpha\int_Mf_iu_ju_j\phi^2(|\nabla u|^2+\epsilon)^{\alpha-1}u_{ki}u_k
\\&-&\int_M\nabla|\nabla u|^2\nabla\phi^2(|\nabla u|^2+\epsilon)^\alpha-\alpha\int_M\nabla|\nabla u|^2\nabla|\nabla u|^2\phi^2(|\nabla u|^2+\epsilon)^{\alpha-1}
\\&=&\int_Mf_iu_ju_j(\phi^2)_i(|\nabla u|^2+\epsilon)^\alpha+2\alpha\int_Mf_iu_ju_j\phi^2(|\nabla u|^2+\epsilon)^{\alpha-1}u_{ki}u_k
\\&-&4\int_M\nabla|\nabla u||\nabla u|\nabla\phi\phi(|\nabla u|^2+\epsilon)^\alpha-4\alpha\int_M|\nabla|\nabla u||^2|\nabla u|^2\phi^2(|\nabla u|^2+\epsilon)^{\alpha-1}
\\&\leq&\int_Mf_iu_ju_j(\phi^2)_i(|\nabla u|^2+\epsilon)^\alpha+2\alpha\int_Mf_iu_ju_j\phi^2(|\nabla u|^2+\epsilon)^{\alpha-1}u_{ki}u_k
\\&-&4\int_M\nabla|\nabla u||\nabla u|\nabla\phi\phi(|\nabla u|^2+\epsilon)^\alpha-4\alpha\int_M|\nabla|\nabla u||^2\phi^2(|\nabla u|^2+\epsilon)^\alpha.
\end{eqnarray*}
Let $\alpha=\frac{p-2}{2}(<0)$. Then the above inequality implies that
\begin{eqnarray*}
&&2(p-1)\int_M|\nabla|\nabla u||^2\phi^2(|\nabla u|^2+\epsilon)^\alpha+\int_MRu_ju_j\phi^2(|\nabla u|^2+\epsilon)^\alpha\nonumber
\\&\leq&\int_Mf_iu_ju_j(\phi^2)_i(|\nabla u|^2+\epsilon)^\alpha+2\alpha\int_Mf_iu_ju_j\phi^2(|\nabla u|^2+\epsilon)^{\alpha-1}u_{ki}u_k
\\&-&4\int_M\nabla|\nabla u||\nabla u|\nabla\phi\phi(|\nabla u|^2+\epsilon)^\alpha
\\&\leq&\int_Mf_iu_ju_j(\phi^2)_i(|\nabla u|^2+\epsilon)^\alpha+2\alpha\int_Mf_iu_ju_j\phi^2(|\nabla u|^2+\epsilon)^{\alpha-1}u_{ki}u_k
\\&-&(p-1)\int_M|\nabla|\nabla u||^2\phi^2(|\nabla u|^2+\epsilon)^\alpha-\frac{4}{p-1}\int_M|\nabla u|^2|\nabla\phi|^2(|\nabla u|^2+\epsilon)^\alpha.
\end{eqnarray*}
Therefore
\begin{eqnarray}\label{ine5}
&&(p-1)\int_M|\nabla|\nabla u||^2\phi^2(|\nabla u|^2+\epsilon)^\alpha+\int_MRu_ju_j\phi^2(|\nabla u|^2+\epsilon)^\alpha\nonumber
\\&\leq&\int_Mf_iu_ju_j(\phi^2)_i(|\nabla u|^2+\epsilon)^\alpha+2\alpha\int_Mf_iu_ju_j\phi^2(|\nabla u|^2+\epsilon)^{\alpha-1}u_{ki}u_k\nonumber
\\&-&\frac{4}{p-1}\int_M\nabla|\nabla u||\nabla u|\nabla\phi\phi(|\nabla u|^2+\epsilon)^\alpha\nonumber
\\&\leq&\int_Mf_iu_ju_j(\phi^2)_i(|\nabla u|^2+\epsilon)^\alpha+2\alpha\int_Mf_iu_ju_j\phi^2(|\nabla u|^2+\epsilon)^{\alpha-1}u_{ki}u_k\nonumber
\\&-&\frac{4}{p-1}\int_M|\nabla u|^2|\nabla\phi|^2(|\nabla u|^2+\epsilon)^\alpha.
\end{eqnarray}
Note that
\begin{eqnarray}\label{eq5}
|\int_Mf_iu_ju_j(\phi^2)_i(|\nabla u|^2+\epsilon)^\alpha|&\leq&2\int_M|\nabla f||\nabla\phi||\nabla u|^p\nonumber
\\&\leq&\frac{2C\sqrt{\lambda}}{r}\int_M|\nabla u|^p\to 0 \ as \ r\to \infty,
\\ \int_M|\nabla u|^2|\nabla\phi|^2(|\nabla u|^2+\epsilon)^\alpha&\leq&\int|\nabla u|^p|\nabla\phi|^2\nonumber
\\ &\leq& \frac{C}{r^2}\int_M|\nabla u|^p \to 0\ as \ r\to \infty.\label{eq6}
\end{eqnarray}
Denote by
$$I=\lim_{r\to\infty}\lim_{\epsilon\to0}\int_Mf_iu_ju_j\phi^2(|\nabla u|^2+\epsilon)^{\alpha-1}u_{ki}u_k.$$
Then from (\ref{eq4}) we get
$$pI=\int_MR|\nabla u|^p\geq0,$$
where we used the fact that $R\geq0$ (see \cite{Chen}). Therefore $I\geq0$, which together with (\ref{ine5})-(\ref{eq6}) imply that
\begin{eqnarray*}
(p-1)\int_M|\nabla|\nabla u||^2|\nabla u|^{p-2}+\int_MR|\nabla u|^p\leq2\alpha I\leq0,
\end{eqnarray*}
where we used $2\alpha=p-2<0$.
Therefore we have $|\nabla u|=c$, which must be zero since $\int_M|\nabla u|^p<\infty$ and $M$ has infinite volume. \endproof

\section{Appendix}
In this appendix we discuss more $L^p$ Liouville type theorems for harmonic or pluriharmonic functions on gradient shrinking or expanding Ricci solitons. To our knowledge, till now no such kind of results were known except for  gradient shrinking K\"ahler-Ricci solitons.

For gradient shrinking Ricci solitons we have
\begin{thm}
Let $(M,g)$ be a gradient shrinking Ricci soliton such that $R\geq\frac{n-2}{2}$. Assume that $u$ is a harmonic function on $M$ with $\int_M|\nabla u|^2<\infty$, then $u$ is a constant function.
\end{thm}
\proof We follow argument used by Munteanu and Sesum in their proof of Theorem 4.1 in \cite{MS}. We assume that $(M,g)$ is complete noncompact. For a cut-off function $\phi$ on $M$ it follows that
\begin{eqnarray*}
\int_MRic(\nabla u,\nabla u)\phi^2&=&-\int_Mf_{ij}u_iu_j\phi^2+\frac12\int_M|\nabla u|^2\phi^2
\\&=&\int_Mu_{ij}f_iu_j\phi^2+\int_Mf_iu_iu_j(\phi^2)_j+\frac12\int_M|\nabla u|^2\phi^2.
\end{eqnarray*}
Note that by integration by parts we have
\begin{eqnarray*}
\int_Mu_{ij}f_iu_j\phi^2=-\frac12\int_M\Delta f|\nabla u|^2\phi^2-\frac12\int_M|\nabla u|^2\langle\nabla f,\nabla\phi^2\rangle.
\end{eqnarray*}
Combing the above two formulas and using $$-\Delta f=R-\frac n2,$$ we get
\begin{eqnarray}\label{app ine1}
\int_MRic(\nabla u,\nabla u)\phi^2&=&\frac12\int_M(R-\frac{n-2}{2})|\nabla u|^2\phi^2-\frac12\int_M|\nabla u|^2\langle\nabla f,\nabla\phi^2\rangle\nonumber
\\&+&\int_M\langle\nabla f,\nabla u\rangle\langle\nabla u,\nabla\phi^2\rangle.
\end{eqnarray}
Recall that we have the following Bochner formula
$$\Delta|\nabla u|^2=2Ric(\nabla u,\nabla u)+2|u_{ij}|^2\geq2Ric(\nabla u,\nabla u)+2|\nabla|\nabla u||^2.$$
Multiplying this inequality by $\phi^2$, using (\ref{app ine1}), and integrating by parts we have
\begin{eqnarray*}
&&2\int_M|\nabla|\nabla u||^2\phi^2+\int_M(R-\frac{n-2}{2})|\nabla u|^2\phi^2
\\&\leq&\int_M|\nabla u|^2\langle\nabla f,\nabla\phi^2\rangle-2\int_M\langle\nabla f,\nabla u\rangle\langle\nabla u,\nabla\phi^2\rangle+\int_M\phi^2\Delta|\nabla u|^2
\\&\leq&\int_M|\nabla u|^2\langle\nabla f,\nabla\phi^2\rangle-2\int_M\langle\nabla f,\nabla u\rangle\langle\nabla u,\nabla\phi^2\rangle
\\&+&\int_M|\nabla|\nabla u||^2\phi^2+4\int_M|\nabla u|^2|\nabla\phi|^2.
\end{eqnarray*}
Therefore we obtain
\begin{eqnarray*}
&&\int_M|\nabla|\nabla u||^2\phi^2+\int_M(R-\frac{n-2}{2})|\nabla u|^2\phi^2
\\&\leq&\int_M|\nabla u|^2\langle\nabla f,\nabla\phi^2\rangle-2\int_M\langle\nabla f,\nabla u\rangle\langle\nabla u,\nabla\phi^2\rangle
+4\int_M|\nabla u|^2|\nabla\phi|^2
\\&\leq&\frac{C}{r^2}\int_M|\nabla u|^2+C\int_{B_p(2r)\setminus B_p(r)}|\nabla u|^2,
\end{eqnarray*}
where in the last inequality we used (\ref{ine00}). Let $r\to \infty$, since $R\geq\frac{n-2}{2}$, we get that $|\nabla u|=c$. We see that $c$ must be zero, i.e. $u$ is a constant function, since $\int_M|\nabla u|^2<\infty$ and $M$ has infinite volume. \endproof

Next we discuss $L^p$ Liouville type theorems for harmonic or pluriharmonic functions on gradient expanding Ricci solitons. Before this, we collect some basic facts of gradient expanding Ricci solitons which will be used in the following. Munteanu and Wang \cite{MW} proved that
$$\inf_{B_p(r)}f\sim-\frac{r^2}{4},$$
as $r(x)$ is large enough, then by $$R+|\nabla f|^2+f=0$$ after a normalization of $f$ (see \cite{ELM}), we have
\begin{eqnarray*}
|\nabla f|^2=-R-f\leq\frac{n}{2}-f\sim\frac{r^2}{4},
\end{eqnarray*}
where we used $R\geq-\frac{n}{2}$ (see \cite{PRS2}\cite{Zha}). Therefore we have
\begin{eqnarray}\label{app ine2}
|\nabla f|\leq\frac{r}{2}+o(r).
\end{eqnarray}
Carrillo and Ni \cite{CN} proved that for any $p\in M$ and $r\geq r_0$, if $R\geq-\sigma$, where $\sigma\geq0$, we have the following monotonicity formula for the volume of gradient expanding Ricci solitons:
\begin{eqnarray}\label{app ine3}
|B_p(r)|\geq|B_p(r_0)|(\frac{r+\tau}{r_0+\tau})^{n-2\sigma},
\end{eqnarray}
where $\tau:=\sqrt{\sigma-f(p)}$.

We have
\begin{thm}\label{app main thm2}
Let $(M,g)$ be a gradient expanding Ricci soliton such that $R\geq-\frac{n-2}{2}$. Assume that $u$ is a harmonic function on $M$ with $\int_M|\nabla u|^2<\infty$, then $u$ is a constant function.
\end{thm}
\proof We follow argument used by Munteanu and Sesum in their proof of Theorem 4.1 in \cite{MS}.  We assume that $(M,g)$ is compete noncompact. For a cut-off function $\phi$ on $M$ it follows that
\begin{eqnarray*}
\int_MRic(\nabla u,\nabla u)\phi^2&=&-\int_Mf_{ij}u_iu_j\phi^2-\frac12\int_M|\nabla u|^2\phi^2
\\&=&\int_Mu_{ij}f_iu_j\phi^2+\int_Mf_iu_iu_j(\phi^2)_j-\frac12\int_M|\nabla u|^2\phi^2.
\end{eqnarray*}
Note that by integration by parts we have
\begin{eqnarray*}
\int_Mu_{ij}f_iu_j\phi^2=-\frac12\int_M\Delta f|\nabla u|^2\phi^2-\frac12\int_M|\nabla u|^2\langle\nabla f,\nabla\phi^2\rangle.
\end{eqnarray*}
Combing the above two formulas and using $-\Delta f=R+\frac n2$, we get
\begin{eqnarray}
\int_MRic(\nabla u,\nabla u)\phi^2&=&\frac12\int_M(R+\frac{n-2}{2})|\nabla u|^2\phi^2-\frac12\int_M|\nabla u|^2\langle\nabla f,\nabla\phi^2\rangle\nonumber
\\&+&\int_M\langle\nabla f,\nabla u\rangle\langle\nabla u,\nabla\phi^2\rangle.
\end{eqnarray}
Recall that we have the following Bochner formula
$$\Delta|\nabla u|^2=2Ric(\nabla u,\nabla u)+2|u_{ij}|^2\geq2Ric(\nabla u,\nabla u)+2|\nabla|\nabla u||^2.$$
Multiplying this inequality by $\phi^2$, using (\ref{app ine1}), and integrating by parts we have
\begin{eqnarray*}
&&2\int_M|\nabla|\nabla u||^2\phi^2+\int_M(R+\frac{n-2}{2})|\nabla u|^2\phi^2
\\&\leq&\int_M|\nabla u|^2\langle\nabla f,\nabla\phi^2\rangle-2\int_M\langle\nabla f,\nabla u\rangle\langle\nabla u,\nabla\phi^2\rangle+\int_M\phi^2\Delta|\nabla u|^2
\\&\leq&\int_M|\nabla u|^2\langle\nabla f,\nabla\phi^2\rangle-2\int_M\langle\nabla f,\nabla u\rangle\langle\nabla u,\nabla\phi^2\rangle
\\&+&\int_M|\nabla|\nabla u||^2\phi^2+4\int_M|\nabla u|^2|\nabla\phi|^2.
\end{eqnarray*}
Therefore we obtain
\begin{eqnarray*}
&&\int_M|\nabla|\nabla u||^2\phi^2+\int_M(R+\frac{n-2}{2})|\nabla u|^2\phi^2
\\&\leq&\int_M|\nabla u|^2\langle\nabla f,\nabla\phi^2\rangle-2\int_M\langle\nabla f,\nabla u\rangle\langle\nabla u,\nabla\phi^2\rangle
+4\int_M|\nabla u|^2|\nabla\phi|^2
\\&\leq&\frac{C}{r^2}\int_M|\nabla u|^2+C\int_{B_p(2r)\setminus B_p(r)}|\nabla u|^2,
\end{eqnarray*}
where in the last inequality we used (\ref{app ine2}). Let $r\to \infty$, since $R\geq-\frac{n-2}{2}$, we get that $|\nabla u|=c$. Since $R\geq-\frac{n-2}{2}$, we see from (\ref{app ine3}) that
 $$|B_p(r)|\geq|B_p(r_0)|(\frac{r+\tau}{r_0+\tau})^2,$$
 which implies that $M$ has infinite volume. Then we have that $c$ must be zero, i.e. $u$ is a constant function, since $\int_M|\nabla u|^2<\infty$. \endproof

Furthermore for gradient expanding K\"ahler-Ricci solitons we have
\begin{thm}
Let $(M,g)$ be a gradient expanding K\"ahler-Ricci soliton. Assume that $u$ is a pluriharmonic function on $M$, $R\geq-\frac{n-2}{2}$, $\int_M|\nabla u|^p<\infty$ for some constant $1<p\leq2$, and $\nabla f(o)=0$ at some point $o\in M$, then $u$ is a constant function.
\end{thm}
\begin{rem}
A point $o\in M$ satisfying $\nabla f(o)=0$ is called an equilibrium point of $(M,g)$, while existence of an equilibrium point of $(M,g)$ plays an important role in rigidity theorems due to Deng and Zhu in \cite{DZ} for complete non-compact gradient expanding or steady K\"ahler-Ricci solitons.
\end{rem}
\proof Since the argument is quite similar with that we have used in the proof of Theorem \ref{main thm steady}, we only give a sketch of proofs here. We only need assume that $1<p<2$. From $u$ is pluriharmonic we see that $\Delta F=0$, where $F=\langle\nabla f, \nabla u\rangle$. Then we can obtain
\begin{eqnarray*}
&&\int_M\phi^2(F^2+\epsilon)^{\frac{p-2}{2}}|\nabla F|^2
\\&\leq&(\frac{2}{p-1})^2\int_M|F|^p|\nabla\phi|^2
\\&\leq&\frac{C}{r^{2-p}}\int_M|\nabla u|^p\nonumber
\end{eqnarray*}
where we used (\ref{app ine2}) to get $$|F|\leq|\nabla f||\nabla u|\leq Cr|\nabla u|.$$ Then by letting $r\to \infty$ we have $\nabla F=0$, which yields that $F$ is a constant function on $M$. Then since $\nabla f(o)=0$ at some point $o\in M$, we see that $F=0$ on $M$.

Let $\alpha<0$ be a constant to be determined later. Then
\begin{eqnarray}\label{eq2'}
&&\int_M Ric(\nabla u,\nabla u)\phi^2(|\nabla u|^2+\epsilon)^\alpha\nonumber
\\&=&\int_Mf_iu_{ij}u_j\phi^2(|\nabla u|^2+\epsilon)^\alpha-\frac12\int_M|\nabla u|^2\phi^2,
\end{eqnarray}
where we used $$f_iu_i=\langle\nabla f,\nabla u \rangle=F=0.$$
Note that
\begin{eqnarray*}
&&\int_Mf_iu_{ij}u_j\phi^2(|\nabla u|^2+\epsilon)^\alpha
\\&=&-\int_M\Delta fu_ju_j\phi^2(|\nabla u|^2+\epsilon)^\alpha-\int_Mf_iu_ju_{ij}\phi^2(|\nabla u|^2+\epsilon)^\alpha
\\&-&\int_Mf_iu_ju_j(\phi^2)_i(|\nabla u|^2+\epsilon)^\alpha-2\alpha\int_Mf_iu_ju_j\phi^2(|\nabla u|^2+\epsilon)^{\alpha-1}u_{ki}u_k.
\end{eqnarray*}
Then by $\Delta f=-R-\frac{n}{2}$ we get
\begin{eqnarray}\label{eq3'}
&&\int_M Ric(\nabla u,\nabla u)\phi^2(|\nabla u|^2+\epsilon)^\alpha\nonumber
\\&=&\frac12\int_M(R+\frac{n}{2}-1)u_ju_j\phi^2(|\nabla u|^2+\epsilon)^\alpha\nonumber
\\&-&\frac12\int_Mf_iu_ju_j(\phi^2)_i(|\nabla u|^2+\epsilon)^\alpha-\alpha\int_Mf_iu_ju_j\phi^2(|\nabla u|^2+\epsilon)^{\alpha-1}u_{ki}u_k.
\end{eqnarray}
Recall that we have the following Bochner formula
\begin{eqnarray*}
&&\Delta|\nabla u|^2
\\&=&2Ric(\nabla u,\nabla u)+2|u_{ij}|^2
\\&\geq&2Ric(\nabla u,\nabla u)+2|\nabla|\nabla u||^2.
\end{eqnarray*}
Then we have
\begin{eqnarray*}\label{ine4'}
&&\int_M\Delta|\nabla u|^2\phi^2(|\nabla u|^2+\epsilon)^\alpha
\\&\geq&2\int_MRic(\nabla u,\nabla u)\phi^2(|\nabla u|^2+\epsilon)^\alpha+2\int_M|\nabla|\nabla u||^2\phi^2(|\nabla u|^2+\epsilon)^\alpha.
\end{eqnarray*}
Plug this inequality into (\ref{eq3'}) we get
\begin{eqnarray*}
&&2\int_M|\nabla|\nabla u||^2\phi^2(|\nabla u|^2+\epsilon)^\alpha+\int_M(R+\frac{n}{2}-1)u_ju_j\phi^2(|\nabla u|^2+\epsilon)^\alpha\nonumber
\\&\leq&\int_Mf_iu_ju_j(\phi^2)_i(|\nabla u|^2+\epsilon)^\alpha+2\alpha\int_Mf_iu_ju_j\phi^2(|\nabla u|^2+\epsilon)^{\alpha-1}u_{ki}u_k
\\&-&4\int_M\nabla|\nabla u||\nabla u|\nabla\phi\phi(|\nabla u|^2+\epsilon)^\alpha-4\alpha\int_M|\nabla|\nabla u||^2\phi^2(|\nabla u|^2+\epsilon)^\alpha.
\end{eqnarray*}
Let $\alpha=\frac{p-2}{2}(<0)$. Then from the above inequality we have
\begin{eqnarray}\label{ine5'}
&&(p-1)\int_M|\nabla|\nabla u||^2\phi^2(|\nabla u|^2+\epsilon)^\alpha+\int_M(R+\frac{n}{2}-1)u_ju_j\phi^2(|\nabla u|^2+\epsilon)^\alpha\nonumber
\\&\leq&\int_Mf_iu_ju_j(\phi^2)_i(|\nabla u|^2+\epsilon)^\alpha+2\alpha\int_Mf_iu_ju_j\phi^2(|\nabla u|^2+\epsilon)^{\alpha-1}u_{ki}u_k\nonumber
\\&-&\frac{4}{p-1}\int_M|\nabla u|^2|\nabla\phi|^2(|\nabla u|^2+\epsilon)^\alpha.
\end{eqnarray}
Note that
\begin{eqnarray}\label{eq5'}
&&|\int_Mf_iu_ju_j(\phi^2)_i(|\nabla u|^2+\epsilon)^\alpha|\nonumber
\\&\leq&2\int_M|\nabla f||\nabla\phi||\nabla u|^p\nonumber
\\&\leq&C\int_{B_p(2r)\setminus B_p(r)}|\nabla u|^p\to 0 \ as \ r\to \infty,
\end{eqnarray}

and
\begin{eqnarray}\label{eq6'}
&&\int_M|\nabla u|^2|\nabla\phi|^2(|\nabla u|^2+\epsilon)^\alpha\nonumber
\\&\leq&\int|\nabla u|^p|\nabla\phi|^2\nonumber
\\ &\leq& \frac{C}{r^2}\int_M|\nabla u|^p \to 0\ as \ r\to \infty.
\end{eqnarray}
Denote by
$$I=\lim_{r\to\infty}\lim_{\epsilon\to0}\int_Mf_iu_ju_j\phi^2(|\nabla u|^2+\epsilon)^{\alpha-1}u_{ki}u_k.$$
Similarly with the proof of Theorem \ref{main thm steady} we get
$$pI=\int_M(R+\frac{n}{2})|\nabla u|^p\geq0,$$
where we used the fact that $R\geq-\frac{n}{2}$. Therefore $I\geq0$, which together with (\ref{ine5'})-(\ref{eq6'}) imply that
\begin{eqnarray*}
(p-1)\int_M|\nabla|\nabla u||^2|\nabla u|^{p-2}+\int_M(R+\frac{n}{2}-1)|\nabla u|^p\leq2\alpha I\leq0,
\end{eqnarray*}
where we used $2\alpha=p-2<0$. Therefore we have $|\nabla u|=c$, which must be zero since $\int_M|\nabla u|^p<\infty$ and $M$ has infinite volume. \endproof

\vspace{0.3cm}

\textbf{Acknowledgments.} The author would like to thank Xueyuan Wan, Jianyu Ou and Yuxing Deng for stimulating  discussions. Many thanks to Shijin Zhang for his interests and useful suggestions which made this paper more readable. Thank you very much the anonymous reviewer for pointing out several crucial mistakes on proofs in the first version of this paper, which force us to add some assumptions in the main theorems. This paper is supported by the Natural Science Foundation of China (Grant no. 12271069).
{}
\vspace{1cm}\sc

Yong Luo

Mathematical Science Research Center of Mathematics,

Chongqing University of Technology,

Chongqing, 400054, China

{\tt yongluo-math@cqut.edu.cn}


\begin{thebibliography}{2}
\bibitem{Cao} H. D. Cao, Recent progress on Ricci solitons, {\em Adv. Lect. Math. (ALM)}, {\bf11},
International Press, Somerville, MA, 2010, 1--38.
\bibitem{CZ} H. D. Cao and D. T. Zhou, On complete gradient shrinking Ricci solitons, {\em J. Differential Geom.} {\bf85} (2010), no. 2, 175--185.
\bibitem{CN} J. A. Carrillo and L. Ni,  Sharp logarithmic Sobolev inequalities on gradient solitons and applications, {\em Comm. Anal. Geom.} {\bf17} (2009), no. 4, 721--753.
\bibitem{Chen} B. L. Chen, Strong uniqueness of the Ricci flow, {\em J. Differential Geom.} {\bf82} (2009), no. 2, 363--382.
\bibitem{ChY} S. Y. Cheng and S. T. Yau,  Differential equations on Riemannian manifolds and their geometric applications, {\em Comm. Pure Appl. Math. } {\bf28} (1975), no. 3, 333--354.
\bibitem{CLN} B. Chow, P. Lu and L. Ni, Hamilton's Ricci flow, {\em Grad. Stud. Math.} {\bf77},
American Mathematical Society, Providence, RIScience Press Beijing, New York, 2006, xxxvi+608 pp.
\bibitem{CM} T. H. Colding and W. P. Minicozzi II, Liouville prepeties, {\em ICCM Not.} {\bf7} (2019), no.1, 16--26.
\bibitem{DZ} Y. X. Deng and X.H. Zhu, Complete non-compact gradient Ricci solitons with nonnegative Ricci curvature, {\em Math. Z.} {\bf279} (2015), no. 1-2, 211--226.
\bibitem{DHT} N. T. Dung, N. T. Le Hai and N. T. Thanh, Eigenfunctions of the weighted Laplacian and a vanishing theorem on gradient steady Ricci soliton, {\em J. Math. Anal. Appl.} {\bf416} (2014), no. 2, 553--562.
\bibitem{ELM} M. Eminenti, G. La Nave and C. Mantegazza,  Ricci solitons: the equation point of view, {\em Manuscripta Math.} {\bf127} (2008), no. 3, 345--367.
\bibitem{GZ} H. B. Ge and S. J. Zhang, Liouville-type theorems on the complete gradient shrinking Ricci solitons, {\em Differential Geom. Appl.} {\bf56} (2018), 42--53.
\bibitem{Ham} R. S. Hamilton, The formation of singularities in the Ricci flow, In:{\em Surveys in Diff. Geom.,} vol.2, International Press, Cambridge, MA, 1995, 7--136.
\bibitem{HM} R. Haslhofer and R. M\"uller, A compactness theorem for complete Ricci shrinkers, {\em Geom. Funct. Anal.} {\bf21} (2011), no. 5, 1091--1116.
\bibitem{Li} P. Li, Harmonic functions and applications to complete manifolds, {\em Lect. notes,} math.uci.edu/pli.
\bibitem{Li2} P. Li, On the structure of complete K\"ahler manifolds with nonnegative curvature near infinity, {\em Invent. Math.} {\bf99} (1990), no. 3, 579--600.
\bibitem{MO} W. X. Mai and J. Y. Ou, Liouville theorem on Ricci shrinkers with constant scalar curvature and its applications, arXiv:2208.07101v1, submitted.
\bibitem{MS} O. Munteanu and N. Sesum,  On gradient Ricci solitons, {\em J. Geom. Anal.} {\bf23} (2013), no. 2, 539--561.
\bibitem{MW} O. Munteanu and J. P. Wang, Analysis of weighted Laplacian and applications to Ricci solitons, {\em Comm. Anal. Geom.} {\bf20} (2012), no.1, 55--94.
\bibitem{MW2} O. Munteanu and J. P. Wang, Holomorphic functions on K\"ahler-Ricci solitons, {\em J. Lond. Math. Soc. (2)} {\bf89} (2014), no. 3, 817--831.
    \bibitem{PRS} S. Pigola, M. Rigoli and A. G. Setti, Vanishing theorems on Riemannian manifolds, and geometric applications, {\em J. Funct. Anal.} {\bf229} (2005), no. 2, 424--461.
    \bibitem{PRS2} S. Pigola, M. Rimoldi and A. G. Setti, Remarks on non-compact gradient Ricci solitons, {\em Math. Z.} {\bf268} (2011), no. 3-4, 777--790.
\bibitem{WW} J. Y. Wu and P. Wu, Harmonic and Schrodinger functions of polynomial growth on gradient shrinking Ricci solitons, {\em Geom. Dedicata} {\bf 217} (2023), no. 4, Paper No. 75, 25 pp.
\bibitem{Yau1} S. T. Yau, Harmonic functions on complete Riemannian manifolds, {\em Comm. Pure Appl. Math.} {\bf28} (1975), 201--228.
\bibitem{Yau2} S. T. Yau, Some function-theoretic properties of complete Riemannian manifold and their applications to geometry, {\em Indiana Univ. Math. J.} {\bf25} (1976), no. 7, 659--670
\bibitem{Zha} S. J. Zhang, On a sharp volume estimate for gradient Ricci solitons with scalar curvature bounded below, {\em Acta Math. Sin. (Engl. Ser.)} {\bf27} (2011), no. 5, 871--882.
\end{thebibliography}
\end{document}